\documentclass[12pt,a4paper,reqno]{amsart}
\usepackage{amsmath}
\allowdisplaybreaks
\usepackage{amsfonts}
\usepackage{amssymb,amsthm,amsfonts,amsthm,latexsym,enumerate,url,cases}
\numberwithin{equation}{section}
\usepackage{mathrsfs}
\usepackage{hyperref}
\hypersetup{colorlinks=true,citecolor=blue,linkcolor=blue,urlcolor=blue}
\newtheorem{theorem*}{Theorem}
\newtheorem{lemma*}{Lemma}
\theoremstyle{plain}

\newtheorem{conjecture}{Conjecture}
\theoremstyle{definition}


\begin{document}

\title
[{A counterexample of two Romanov type conjectures}] {A counterexample of two Romanov type conjectures}

\author
[Yuchen Ding] {Yuchen Ding}

\address{(Yuchen Ding) School of Mathematical Science,  Yangzhou University, Yangzhou 225002, People's Republic of China}
\email{ycding@yzu.edu.cn}


\begin{abstract} In this note, we disprove two Romanov type conjectures posed by Chen.
\end{abstract}
\maketitle

\baselineskip 18pt

\section{Introduction}
For any subset $A$ of natural numbers set $\mathbb{N}$, let $A(x)=|A\cap[1,x]|$.
In 2008, Chen communicated to the authors of \cite{LP} and made the following two conjectures therein. 
\begin{conjecture}\label{A}
Let $\mathscr{A}$ and $\mathscr{B}$ be two sets of positive integers. If there exists a constant $c>0$ such that $\mathscr{A}(\log x/\log 2)\mathscr{B}(x)>cx$ for all sufficiently large $x$, then the set $\{2^a+b:a\in \mathscr{A},~b\in \mathscr{B}\}$ has a positive lower asymptotic density.
\end{conjecture}

\begin{conjecture}\label{B}
Let $\mathscr{A}$ and $\mathscr{B}$ be two sets of positive integers. If there exists a constant $c>0$ such that $\mathscr{A}(\log x/\log 2)\mathscr{B}(x)>cx$ for infinitely many $x$, then the set $\{2^a+b:a\in \mathscr{A},~b\in \mathscr{B}\}$ has a positive upper asymptotic density.
\end{conjecture}

The Ramanov theorem \cite{Ro} offers a positive answer to the conjectures of Chen for $\mathscr{B}$ being the set of primes. However, it will not be the case for general situations.
The disavowals of Conjecture \ref{A} and \ref{B} will be illustrated simultaneously by constructing a concrete counterexample in next section.  

\section{Counterexample}

Let $p_i$ be the $i$--th odd prime and $d_t=p_1p_2\cdots p_t$ for $t\in \mathbb{Z}_+$. For any $t\in \mathbb{Z}_+$, define
$$\mathscr{B}_t=\{n:n\in\mathbb{N},~d_t|n\}\cap\left[2^{2^{t^2}},2^{2^{(t+1)^2}}\right),~~\mathscr{B}=\bigcup_{t=1}^{\infty}\mathscr{B}_t$$
and
$$\mathscr{A}=\mathbb{N},~~\mathscr{C}=2^{\mathscr{A}}+\mathscr{B}=\{2^a+b:a\in \mathscr{A},~b\in \mathscr{B}\}.$$
Let $x$ be a large number and $j$ be the number such that
\begin{align}\label{1}
2^{2^{j^2}}\le x<2^{2^{(j+1)^2}},
\end{align}
which means that
\begin{align}\label{2}
\sqrt{\log\log x}<j\le2\left(\sqrt{\log 2}\right)^{-1}\sqrt{\log\log x}.
\end{align}
From the Chebyshev estimate, we have
\begin{align}\label{3}
d_{j}&=\exp\left(\sum_{2<p\le p_j}\log p\right)\nonumber\\
&\le\exp\left(2j\log j\right)\nonumber\\
&\le\exp\left(3\sqrt{\log\log x}\log\log\log x\right).
\end{align}
In view of equation (\ref{1}), we know $2^{2^{(j-1)^2}}<\sqrt{x}$. So by the construction of $\mathscr{B}$ and equation (\ref{3}), we have 
\begin{align}\label{4}
\mathscr{B}(x)&\ge \frac{x-2^{2^{j^2}}}{d_{j}}+\frac{2^{2^{j^2}}-2^{2^{(j-1)^2}}}{d_{j-1}}-2\nonumber\\
&\ge\frac{x-2^{2^{(j-1)^2}}}{d_{j}}-2\nonumber\\
&\gg \frac{x}{\exp\left(3\sqrt{\log\log x}\log\log\log x\right)}.
\end{align}
It follows that
\begin{align}\label{5}
\mathscr{A}(\log x/\log 2)\mathscr{B}(x)\gg \frac{x\log x}{\exp\left(3\sqrt{\log\log x}\log\log\log x\right)}>x.
\end{align} 
It remains to prove $\mathscr{C}(x)=o(x)$ as $x\to\infty$. It is clear that
\begin{align}\label{6}
\mathscr{C}(x)=\#\left\{c\le x:c=2^a+b,a\in \mathscr{A},b\in\mathscr{B}\right\}\le S_1(x)+S_2(x),
\end{align}
where
\begin{align}\label{7}
S_1(x)=\#\left\{c\le x:c=2^a+b,a\in \mathscr{A},b\in\mathscr{B}_{j}\right\}
\end{align}
and
\begin{align}\label{8}
S_2(x)=\#\left\{c\le x:c=2^a+b,a\in \mathscr{A},b\notin\mathscr{B}_{j}\right\}.
\end{align}
Note that if $c=2^a+b$ for some $b\in\mathscr{B}_{j}$, then $p_i\nmid c$ for any $1\le i\le j$. This fact leads to the following bound
\begin{align}\label{9}
S_1(x)&\le\sum_{\substack{c\le x\\\left(c,\prod_{p\le p_j}p\right)=1}}1
\nonumber\\
&=\sum_{\ell|\prod_{p\le p_j}}\mu(\ell)\left\lfloor\frac{x}{\ell}\right\rfloor\nonumber\\
&\le x\prod_{p\le p_j}\left(1-\frac{1}{p}\right)+2^j.
\end{align}
The same observation yields the following estimate
\begin{align}\label{10}
S_2(x)&\le \#\left\{c\le x:c=2^a+b,a\in \mathscr{A},b\in\mathscr{B}_{j-1}\right\}+2^{2^{(j-1)^2}}\frac{\log x}{\log 2}\nonumber\\
&\ll x\prod_{p\le p_{j-1}}\left(1-\frac{1}{p}\right)+\sqrt{x}\log x+2^j.
\end{align}
It can be seen that 
\begin{align}
2^j\ll \exp\left(3\sqrt{\log\log x}\right)\ll\sqrt{x}
\end{align}
 from equation (\ref{2}). 
Therefore, combing equations (\ref{6}), (\ref{9}) and (\ref{10}) we have
\begin{align}\label{11}
\mathscr{C}(x)&\ll x\prod_{p\le p_{j-1}}\left(1-\frac{1}{p}\right)+\sqrt{x}\log x\nonumber\\
&\ll x\left(\log p_{j-1}\right)^{-1}+\sqrt{x}\log x\nonumber\\
&\ll x\left(\log\log\log x\right)^{-1},
\end{align}
where the last but one step follows from the Mertens estimate, which is  surely more to expectation than our requirement.

\section{Remarks}
The Romanov type problems start from the remarkable paper \cite{Ro} of Romanv, where he proved that there is a positive lower asymptotic density of odd numbers which can be represented by the sum of a prime and a power of $2$. 

On the opposite direction, van der Corput \cite{va} showed that there is a positive lower asymptotic density of odd numbers none of whose members can be represented by the sum of a prime and a power of $2$. Soon, Erd\H os \cite{Er} constructed an arithmetic progression of odd numbers having the same property required in the paper of van der Corput. The results of van der Corput and Erd\H os give a negative answer to an old conjecture of de Polignac \cite{de1,de2}.

At present, the author of this note has no answer to the following question. Let $\mathscr{B}$ be a set of positive integers satisfying $\mathscr{B}(x)=O(x/\log x)$, then is it true that
$$\limsup_{x\to \infty}\frac{\mathscr{C}(x)}{\mathscr{B}(x)}=\infty,$$
where $$\mathscr{C}=2^{\mathbb{N}}+\mathscr{B}=\{2^k+b:k\in \mathbb{N},~b\in \mathscr{B}\}.$$

\section*{Acknowledgments}
The author is
supported by the Natural Science Foundation of Jiangsu Province of
China, Grant No. BK20210784. He is also supported by the
foundations of the projects ``Jiangsu Provincial
Double--Innovation Doctor Program'', Grant No. JSSCBS20211023 and
``Golden  Phoenix of the Green City--Yang Zhou'' to excellent PhD,
Grant No. YZLYJF2020PHD051.

\end{document}